\documentclass[10pt,twoside]{article}
\usepackage{amsmath,amssymb,graphicx,latexsym,array}

\newtheorem{thm}{Theorem}
\newtheorem{rem}{Remark}
\newtheorem{exa}{Example}

\newcommand{\diag}{\operatorname{diag}}
\newcommand{\sgn}{\operatorname{sgn}}
\newcommand{\co}{\operatorname{\overline{co}}}
\newcommand{\cappa}{\operatorname{K}}
\newcommand{\ds}{\displaystyle}

\def\qed{\hfill $\square$}
\def\squar{\vbox{\hrule\hbox{\vrule height 6pt \hskip
6pt\vrule}\hrule}}

\def\qed{\hfill $\squar$}
\def\squar{\vbox{\hrule\hbox{\vrule height 6pt \hskip
6pt\vrule}\hrule}}

\makeatletter \oddsidemargin.9375in \evensidemargin \oddsidemargin
\begin{document}
\thispagestyle{empty} \setcounter{page}{1}

\noindent
{\footnotesize {\bf To appear in\\[-1.75mm]
{\em Dynamics of Continuous, Discrete and Impulsive Systems}}\\[-1.00mm]
http:monotone.uwaterloo.ca/$\sim$journal} $~$ \\ [.3in]


\begin{center}
{\large\bf Synchronization problems for unidirectional \\feedback
coupled nonlinear systems} \vskip.20in

{\large  Oleg Makarenkov$^{1}$, Paolo Nistri$^{2}$,
Duccio Papini$^{2}$}\\[2mm]

{\footnotesize
$^{1}$ Dept. of Mathematics, Voronezh State University, Voronezh, Russia\\
e-mail: omakarenkov@kma.vsu.ru\\

$^{2}$ Dip. di Ingegneria dell' Informazione, Universit\`a di Siena, 53100 Siena, Italy\\
e-mail: \{pnistri, papini\}@dii.unisi.it }
\end{center}

{\footnotesize \noindent {\bf Abstract.}  In this paper we
consider three different synchronization problems consisting in
designing a nonlinear feedback unidirectional coupling term for
two (possibly chaotic) dynamical systems in order to drive the
trajectories of one of them, the slave system, to a reference
trajectory or to a prescribed neighborhood of the reference
trajectory of the second dynamical system: the master system. If
the slave system is chaotic then synchronization can be viewed as
the control of chaos; namely the coupling term allows to suppress
the chaotic motion by driving the chaotic system to a prescribed
reference trajectory. Assuming that the entire vector field
representing the velocity of the state can be modified, three
different methods to define the nonlinear feedback synchronizing
controller are proposed: one for each of the treated problems.
These methods are based on results from the small parameter
perturbation theory of autonomous systems having a limit cycle,
from nonsmooth analysis and from the singular perturbation theory
respectively. Simulations to illustrate the effectiveness of the
obtained results are also presented.}








\vskip.2in

\section{Introduction}

\noindent In recent years in the literature on dynamical system
analysis a considerable attention has been devoted to the problem
of synchronization of coupled nonlinear dynamical systems (see
e.g. \cite{ACH}, \cite {ble}, \cite{La1}, \cite{PRK}, \cite{Stra},
\cite{WW}). One of the most effective methods for solving such
problem consists in designing a feedback coupling term which
drives the trajectories of one of the two systems (the so-called
slave system) to a prescribed reference trajectory of the second
one (named master system). Examples of such approach can be found,
for instance, in \cite{STDMG1} where the coupling term is
represented by a linear feedback of the tracking error. In
\cite{YZ} a bidirectional linear coupling term is proposed to
synchronize two chaotic systems. An approach to synchronization
based on the classical notion of observers, when the state is not
fully available, can be found in \cite{feki1} and \cite{nm}. In
many cases, when one deals with nonlinear chaotic dynamical
systems, the interest is that of steering any trajectory of the
chaotic system to an equilibrium point or to a limit cycle of the
same system or of another coupled system, see \cite{CD},
\cite{FP}, \cite{LM} and \cite{WC1}. For an adaptive control
approach using a linear reference model we refer to \cite{W}.
Finally, in \cite{La2} it is stressed how the carelessness
application of the mathematical tools of the synchronization
theory can lead to incorrect results.

In this paper, under the condition that the entire vector field of
the velocity can be modified, we aim at designing a nonlinear
feedback unidirectional coupling term, based on the state model
system, in such a way that all the trajectories of the slave
system are steered to a prescribed reference trajectory of the
master system. In other words, the coupling term makes stable, in
a sense that it will be precised for each problem in the
following, a prescribed trajectory of a dynamical system with
respect to the trajectories of an other dynamical system by
coupling these systems by means of a suitably defined nonlinear
feedback coupling term.

Following the linear feedback approach of the previously cited
references, we provide examples of how it can be possible, by
means of different mathematical theories, to define a nonlinear
unidirectional feedback coupling term in order to determine a
prescribed dynamical behavior to some classes of nonlinear
dynamical systems in the case when this term can affect each
component of the velocity vector field.

The problems that we will treat in this paper are illustrated in
the sequel. The first one is the problem of the synchronization of
the phase of a limit cycle of an autonomous system with that of
the limit cycle of the same period of another autonomous system.
The feedback design is based on classical results due to Malkin
\cite{mal2} on the existence of periodic solutions of an
autonomous system perturbed by a small parameter nonautonomous
term and on their behavior when the perturbation disappears,
namely when the parameter tends to zero. Many authors, see for
instance \cite{bokbk}, \cite{iosk} and \cite{rpk} and the
extensive references therein, have considered the problem of the
control of the balance between the phases of the subsystems state
variables oscillations by coupling the subsystems in different
ways, i.e. by suitably balancing the energy due to the
interaction. In particular in \cite{bokbk} a dynamic feedback
coupling term for phase locking of non identical oscillators is
presented.

The second and third problem consist in the synchronization of the
trajectories to a reference trajectory of the master. To solve
these problems we adopt two different feedback laws based on a
sliding manifold approach. First, we define a static discontinuous
feedback coupling term with a gain depending on the bounded set of
the initial conditions for the trajectories of the slave system.
It is defined by means of the signum function of the tracking
error, that is by the signum of the difference between a
trajectory of the slave system and the reference trajectory. By
means of a suitably defined nondifferentiable Liapunov function
and its subdifferentiability properties \cite{clarke} we can prove
that any trajectory originating from a given bounded set converges
to the reference trajectory in an estimated finite time which
depends on the feedback gain. It is worth to observe that, since
the right hand side of the slave system is discontinuous with
respect to the state, it is necessary to introduce a suitable
concept of solution for this system, in fact we consider solutions
in the sense of Filippov \cite{filippov1}. The discontinuity along
the reference curve $y_0=y_0(t), t\ge 0,$ of the feedback law
makes the coupled system robust against modelling errors and
external disturbances, in the sense that the tracking error tends
to zero in finite time also in presence of modelling imprecision
and disturbances, if the gain is sufficiently large. Moreover,
observe that it is possible to get the reference trajectory at any
prescribed speed by suitably increasing the gain. Since the
implementation of the associated switchings across the reference
curve is necessarily imperfect, in practice switching is not
instantaneous and the value of the tracking error
$e(t)=x(t)-y_0(t), t\ge 0,$ is not perfectly known; this leads to
the chattering phenomenon which is the main drawback of this
feedback law, namely the trajectory of the slave system rapidly
oscillates around the reference trajectory. This is a quite
undesirable effect, in fact each component of the signum of the
tracking error in the coupling term switches very fast between the
positive and negative value of the corresponding gain, which is
not a feasible behavior for the physical implementation of the
control law. In the framework of synchronization of chaotic system
a feedback of this type has been used in \cite{WC1}, where the
chattering phenomenon has been also emphasized.

To avoid the undesirable chattering phenomenon we then consider a
dynamic feedback, as introduced in \cite{CDN1}, defined by means
of a differential equation involving the slave system and the
reference trajectory. This equation depends also on a small
parameter and it satisfies, under general conditions, all the
assumptions of the classical singular perturbation theory on
infinite intervals \cite{hop1}. As we will see this ensures that
any trajectory of the slave system approaches the reference
trajectory within any prescribed error.

Finally, we present a simulation for each of the considered
problems which illustrates the effectiveness of the obtained
results. Precisely, for the first problem we consider, both for
the master and slave system, the same FitzHugh-Nagumo type
equation which has an asymptotically stable limit cycle and by
implementing our method we synchronize the phase of the slave with
that of the master. For the second and third problem we have
considered a chaotic neural network as the slave system
\cite{CZL}, and as master system a neural network which possesses
a globally asymptotically exponentially stable periodic solution
which is taken as reference trajectory \cite{PT}.

The paper is organized as follows. In Section \ref{sec:phase} we
treat the problem of the phase synchronization of two
self-oscillating nonlinear dynamical systems. In Section
\ref{sec:static} for two nonautonomous systems we steer any
trajectory of one of these two systems to a prescribed trajectory
of the second one by means of a high gain discontinuous feedback
of the tracking error. In Section \ref{sec:dynamic} we consider
the same systems and we design a dynamical feedback coupling term
which drives any trajectory of the slave system to any prescribed
neighborhood of the reference trajectory. Finally, in Section
\ref{sec:examples} we present the simulations which illustrate the
meaning of the obtained results.

\section{Phase synchronization of limit cycles}\label{sec:phase}

\noindent
In this section we consider two autonomous systems
\begin{equation}\label{eq1}
  \dot x = f(x),
\end{equation}
where $f\in C^2(\mathbb{R}^n,\mathbb{R}^n),$ and
\begin{equation}\label{eq2}
  \dot y = g(y),
\end{equation}
where $g\in C^2(\mathbb{R}^n,\mathbb{R}^n).$ We assume that they
have orbitally stable limit cycles $x_0(t)$ and $y_0(t),$
respectively, of the same period $T.$ According to the terminology
adopted in the literature for the problem that we will treat here,
we will refer to systems (\ref{eq1}) and (\ref{eq2}) as slave and
master system respectively. We are interested in the phase
synchronization of the slave system to that of the master system
by adding a coupling term to (\ref{eq1}). Specifically, we will
show that for every $\mu>0$ we can define a function
$\Phi_\mu:\mathbb{R}^n\times\mathbb{R}^n\to\mathbb{R}^n$ such that
the system
\begin{equation}\nonumber
  \dot x = f(x)+\Phi_\mu(x,y_0(t)),\quad t\ge 0,
\end{equation}
has a unique asymptotically stable $T$-periodic solution $x_\mu$
in a neighborhood of $\{x_0(t):t\in[0,T]\}$ satisfying the property
\begin{equation}\label{prop}
\left|\int_0^T\left| x_\mu(\tau)-y_0(\tau)\right|^2
d\tau-\min_{s\in[0,T]}\int_0^T\left|
x_0(\tau+s)-y_0(\tau)\right|^2 d\tau\right|<\mu.
\end{equation}
For this we assume that the Floquet multiplier of the linearized
systems around $x_0(t)$ and $y_0(t),$ equal to 1 is simple and
that the others $n-1$ are inside of the unit open circle. We
consider the following system
\begin{eqnarray}\label{peq}
\begin{aligned}
   \dot x & =
f(x)+\varepsilon\Big(|x-y_0(t)|^2-\\
&-\min_{s\in[0,T]}\frac{1}{T}\int_0^T | x_0(\tau+s)-y_0(\tau)|^2
d\tau-\delta\Big)f(x),
\end{aligned}
\end{eqnarray}
where $\varepsilon$ and $\delta$ are positive scalar parameters.

\vskip0.1truecm\noindent We can prove the following result.

\begin{thm}\label{thm:phase}
Assume, that the equation
\begin{equation}\label{cond}
\int_0^T | x_0(\tau+\theta)-y_0(\tau)|^2
d\tau=\min_{s\in[0,T]}\int_0^T| x_0(\tau+s)-y_0(\tau)|^2 d\tau
\end{equation} has an unique solution $\theta_0\in[0,T].$ Then for every $\mu>0$ there exists
$\delta_\mu>0$ such that for every $\delta\in[0,\delta_\mu]$ there
is $\varepsilon_\delta>0$ for which the following results hold for
$\varepsilon\in(0,\varepsilon_\delta).$

\noindent 1) System (\ref{peq}) possesses a unique asymptotically
stable $T$-periodic solution $x_\mu$ such that
$$  x_\mu(t)\in
  \mathcal{N}_{\delta_\mu}(x_0),\ {\rm for\ any\ }t\in[0,T],
$$
where
$\mathcal{N}_{\delta_\mu}(x_0)=\left\{x\in\mathbb{R}^n:\ds\inf_{t\in[0,T]}|x-x_0(t)|<\delta_\mu\right\}$
denotes the $\delta_\mu$-neighbor\-hood in $\mathbb{R}^n$ of the
limit cycle $x_0,$

\noindent 2) the solution $x_\mu$ satisfies property (\ref{prop}).
\end{thm}

\noindent To prove this theorem we need the following result due
to I.~G.~Malkin \cite{mal2}, which is one of the main tools for
the study of the synchronization of coupled systems (see
\cite{ble}).  Consider the system
\begin{equation}\label{mal1}
  \dot x  = f(x)+\varepsilon \gamma(t,x)
\end{equation}
where $\gamma\in{C^1(\mathbb{R}\times
\mathbb{R}^n,\mathbb{R}^n)},$ assume that $\gamma$ is $T$-periodic
with respect to time. Then it is possible to show that system
$$
  \dot z=-\left(f'(x_0(t))\right)^*z
$$
has a $T$-periodic solution $z^*$ such that
\begin{equation}\label{good}
  \left<z^*(t),\dot x_0(t)\right>=1,\quad {\rm for\ any\ } t\in[0,T].
\end{equation}
Let us introduce the function $F:\mathbb{R}\to \mathbb{R}$ as
follows
$$
  F(\theta)=\int_0^T
  \left<z^*(\tau),\gamma(\tau-\theta,x_0(\tau))\right>d\tau,\quad
  {\rm for\ any\ }\theta\in\mathbb{R}.
$$
We can now formulate the following result.
\begin{thm}\label{thm:mal}{\rm \bf(\cite{mal2}, Theorems pp. 387 and 392)}
Assume that for sufficiently small $\varepsilon>0$ system
(\ref{mal1}) has a continuous family $\varepsilon\to
x_\varepsilon$ of $T$-periodic solutions satisfying the property
\begin{equation}\label{rel}
  x_\varepsilon(t)\to x_0(t+\theta_0){\ \ \rm as\ \ }\varepsilon\to 0
\end{equation}
then $F(\theta_0)=0.$ Moreover, if $F(\theta_0)=0$ and
$F'(\theta_0)\not=0$ then (\ref{rel}) holds true and the solutions
$x_\varepsilon$ are asymptotically stable or unstable according to
whether $F'(\theta_0)$ is negative or positive.
\end{thm}

\noindent {\bf Proof of Theorem 1.} For $\delta>0$ let
$$
\gamma_\delta(t,x)=\left(|x-y_0(t)|^2-\min_{s\in[0,T]}\frac{1}{T}\int_0^T|
x_0(\tau+s)-y_0(\tau)|^2 d\tau-\delta\right)f(x).
$$
Observe that $\gamma_\delta\in C^1([0,T]\times
\mathbb{R}^n,\mathbb{R}^n)$ is $T$-periodic with respect to time
and so it can be extended from $[0,T]$ to $\mathbb{R}$ by
$T$-periodicity. By  (\ref{good}) we have now that
\begin{eqnarray}
\begin{aligned}
  F_\delta(\theta) & = \int_0^T |x_0(\tau+\theta)-y_0(\tau)|^2
  d\tau-\\
  & -\min_{s\in[0,T]}\int_0^T |x_0(\tau+s)-y_0(\tau)|^2 d\tau-T\delta.
\end{aligned}
\end{eqnarray}
Under our assumptions the function $\theta\to F_\delta(\theta)$ is
$T$-periodic and continuously differentiable. Furthermore there
exists a unique $\theta_0\in[0,T]$ such that $F_0(\theta_0)=0.$
Since the function $F_0$ reaches its minimum at $\theta_0$, we
have that $F_0'(\theta_0)=0.$ Then, there exists $\nu_0$ such that
$F_0'(\theta)\not=0$ for any
$\theta\in(\theta_0-\nu_0,\theta_0)\cup(\theta_0,\theta_0+\nu_0).$
Thus  $F_0'(\theta)<0$ for $\theta\in(\theta_0-\nu_0,\theta_0)$
and $F_0'(\theta)>0$ for $\theta\in(\theta_0,\theta_0+\nu_0).$
Therefore for any $\delta>0$ sufficiently small there exists a
unique $\theta_\delta$ such that

\noindent 1) $F_\delta(\theta_\delta)=0$ and
$F_\delta'(\theta_\delta)>0,$

\noindent 2) $\theta_\delta\to \theta_0$ as $\delta\to 0.$

\noindent Let $\mu>0,$ by the previous considerations there exists
$\delta_\mu>0$ such that
\begin{equation}\label{rela1}
\left|\int_0^T|x_0(\tau+\theta_\delta)-y_0(\tau)|^2 d\tau-
\int_0^T|x_0(\tau+\theta_0)-y_0(\tau)|^2
d\tau\right|<\frac{\mu}{2}
\end{equation}
for any $\delta\in[0,\delta_\mu].$ Then for a given
$\delta\in[0,\delta_\mu],$ by applying the Malkin's results
(Theorem 2 above) and by taking into account (\ref{rel}), we
deduce the existence of a positive number $\varepsilon_\delta>0$
such that, for any $\varepsilon\in[0,\varepsilon_\delta],$ system
(\ref{peq}) possesses a unique asymptotically stable $T$-periodic
solution $x_\mu$ such that $x_\mu(t)\in
\mathcal{N}_{\delta_\mu}(x_0),$
 for any $t\in [0,T],$ satisfying
\begin{equation}\label{rela2}
\left|\int_0^T|x_\mu(\tau)-y_0(\tau)|^2 d\tau-
\int_0^T|x_0(\tau+\theta_\delta)-y_0(\tau)|^2
d\tau\right|<\frac{\mu}{2}.
\end{equation}
Finally, (\ref{rela1}) and (\ref{rela2}) conclude the proof. \qed

\begin{rem}
In the case when equation (\ref{cond}) has $k$ solutions on the
interval $[0,T]$ we can reformulate Theorem 1 by simply replacing
the sentence ``unique asymptotically stable $T$-periodic
solution'' by ``$k$ asymptotically stable $T$-periodic
solutions''. Furthermore, we can interchange the r\^ole of the
slave and master system in all of the previous arguments.
\end{rem}

\section{A static feedback for the synchronization of trajectories}\label{sec:static}

\noindent
In this section we consider two nonlinear nonautonomous
systems
\begin{eqnarray}
& \dot x = \phi(t,x) \label{1} \\
& \dot y = \psi(t,y)\label{2},
\end{eqnarray}
where $\phi,\psi:[0,+\infty)\times \mathbb{R}^n\to\mathbb{R}^n$
satisfy the following conditions:
\begin{itemize}
\item[(3.1)] $t \to\phi(t,x)$ and $t\to\psi(t,x)$ are Lebesgue
measurable functions for any $x\in\mathbb{R}^n;$ $x\to\phi(t,x)$
and $x\to\psi(t,x)$ are locally Lipschitz functions for almost all
(a.a.) $t\ge 0.$
\item[(3.2)] For any $\rho>0$ there exists an
integrable function $\gamma_\rho(t),$ $t\ge 0,$ such that
$$
  |\phi(t,x)|\le\gamma_\rho(t)\quad {\rm and}\quad
  |\psi(t,x)|\le\gamma_\rho(t)
$$
for a.a. $t\ge 0$ and any $x\in\mathbb{R}^n$ such that $|x|\le
\rho.$
\item[(3.3)] Any local solution $x=x(t)$ and $y=y(t)$ of
the Cauchy problems
$$
\left\{\begin{array}{l} \dot x=\phi(t,x)\\
x(0)=x_0\in\mathbb{R}^n
\end{array}\right. \quad
\left\{\begin{array}{l} \dot y=\psi(t,y)\\
y(0)=y_0\in\mathbb{R}^n,
\end{array}\right.
$$
can be extended to the interval $[0,+\infty).$
\end{itemize}
Let $y_0(t), t\ge 0,$ be the prescribed bounded solution of
(\ref{2}) to which system (\ref{1}) must be synchronized by means
of a suitable feedback coupling term to be added to (\ref{1}). The
coupling term we propose here has the form of a static feedback
given by
\begin{equation}
B\sgn(x-y_0(t)) \label{3},
\end{equation}
where $B=\diag(b_i),$ $b_i<0$ for any $i=1,2,...,n,$ and
$\sgn(x)=(\sgn(x_1),\dots,\\ \sgn(x_n))$ for all
$x=(x_i)_{i=1}^n\in\mathbb{R}^n.$ Therefore system \eqref{1} takes
the form
\begin{equation}\label{1p}
\dot x=\phi(t,x)+B\sgn(x-y_0(t)).
\end{equation}
Since the right hand side of system \eqref{1p} is discontinuous in
the state variable $x$ we must adopt a suitable notion of solution
for any Cauchy problem associated to \eqref{1p}. Here, following
\cite{filippov1}, as a solution of \eqref{1p} we intend an
absolutely continuous function $x(t),$ $t\ge 0,$ such that
$$
  \dot x(t)\in\phi(t,x(t))+B\cappa[\sgn(x(t)-y_0(t))]
$$
for a.a. $t\ge 0,$ where
$$
\cappa[\sgn(x-y)]=\bigcap_{\varepsilon>0}\bigcap_{\mu(N)=0}\co
\left(\sgn(\mathcal{B}_{\varepsilon}(x-y)\setminus N)\right),
$$
$N\subset \mathbb{R}^n$ is an arbitrary set of Lebesgue measure
zero, $\co(A)$ denotes the closure of the convex hull of the set
$A$ and $\mathcal{B}_{\varepsilon}(c)$ is the open ball in
$\mathbb{R}^n$ with radius $\varepsilon$ and center $c.$ It is
immediate to see that if $t_0$ is a discontinuity point for some
component of the vector $\sgn(x(t)-y_0(t)),$ then the
corresponding component of $\cappa[\sgn(x(t_0)-y_0(t_0))]$ is the
interval $[-1,1]$. The other component of
$\cappa[\sgn(x(t_0)-y_0(t_0))]$ being $+1$ or $-1.$ In the sequel
we will refer to such solution as a Filippov solution to system
\eqref{1p}.

\noindent Then we associate  to system \eqref{1p} the Lyapunov
function $V:\mathbb{R}^n\to \mathbb{R}_+:=\{x\in \mathbb{R}:
x\ge0\}$ defined as follows
$$
 V(e)=\sum_{i=1}^n|e_i|,
$$
where $e=(e_i)^n_{i=1}$ is the tracking error given by
$e_i=x_i-y_i,$ $i=1,2,...,n.$

\vskip0.1truecm\noindent We can prove the following result.
\begin{thm}\label{thm:static}
For any bounded set $I\subset \mathbb{R}^n$ there exists a
positive constant $M_I>0$ such that, if $b_i<-M_I$ for any
$i=1,2,...,n,$ then any Filippov solution $x(t),$ $t\ge 0,$ to
system \eqref{1p} starting from a point $x_0\in I$ converges in
finite time $t_0\le-\dfrac{V(e(0))}{\mu_I},$ where $\mu_I={\rm
max}\{M_I+b_i,\ i=1,2,...,n \}<0,$ to $y_0(t)$ and $x(t)=y_0(t)$
for any $t\ge t_0.$
\end{thm}

\noindent {\bf Proof.} Let $x(t),$ $t\ge 0,$ be a solution to
\eqref{1p} with $x(0)=x_0,$ whenever $x_0\in I.$ Let
$e(t)=\left(e_i(t)\right)^n_{i=1},$ $e_i(t)=x_i(t)-y_{0,i}(t),$
$t\ge 0.$ From the chain rule for locally Lipschitz regular maps
(\cite{clarke}, Theorem 2.3.9-(iii)), for a.a. $t\ge 0$, we have
that
$$
\frac{d}{dt} V(e(t))=\left<\xi,\dot e(t)\right>,
$$
for any $\xi\in\partial V(e(t)),$ where $\partial V(e(t))$ is the
Clarke generalized gradient of the function $V$ evaluated at
$e(t)$ (\cite{clarke}, p. 27). In this case it can be easily seen
that the vector $\xi=(\xi_i)^n_{i=1}$ is given by $\xi_i=+1,$
$(\xi_i=-1),$ for those indexes $i\in J_+(t), (J_-(t)),$ for which
$e_i(t)>0,$ ($e_i(t)<0$), while $\xi_i\in[-1,1]$ for the indexes
$i\in J_0(t)$ for which $e_i(t)=0.$ Therefore, for a.a. $t\ge 0$
such that $e(t)\not=0,$ we have that
\begin{equation}\label{4_}
\begin{array}{rcl}
\begin{aligned}
 & \dfrac{d}{dt}V(e(t)) = \ds\sum_{i\in J_+(t)\cup J_-(t)}\dot e_i(t) \sgn e_i(t)\\
&= \ds\sum_{i\in J_+(t)\cup J_-(t)}\left[\phi_i(t,x(t))-\psi_i(t,y_0(t))+b_i\sgn e_i(t)\right]\sgn e_i(t)\\
&= \ds\sum_{i\in J_+(t)\cup J_-(t)}\left[\phi_i(t,x(t))-\psi_i(t,y_0(t))\right]\sgn e_i(t)+b_i\\
&\le \max\{M_I+b_i,\ i=1,2,...,n\}=\mu_I, \end{aligned}
\end{array}
\end{equation}
where $M_I\ge|\phi_i(t,x(t))-\psi_i(t,y_0(t))|$ for a.a. $t\ge 0$
and for any $i=1,2,..,n.$ Observe that, under our assumptions on
the vector fields $\phi(t,x)$ and $\psi(t,x),$ such a constant
$M_I>0$ does exist. In fact, $y_0(t),$ $t\ge 0,$ is a bounded
trajectory in $\mathbb{R}^n$ and it is not hard to see that any
solution $x(t),$ $t\ge 0,$ to (\ref{1p}) starting from the bounded
set $I\subset {\mathbb R}^n$ is also bounded in $\mathbb{R}^n.$ By
our assumption $\mu_I<0,$ therefore by integrating (\ref{4_})
between $0$ and $t>0$ we obtain
$$
V(e(t))-V(e(0))\le \mu_I t.
$$
In conclusion, for $t\ge-\dfrac{V(e(0))}{\mu_I}>0$ we have that
$V(e(t))=0$ and thus $e(t)=0$ for any
$t\ge-\dfrac{V(e(0))}{\mu_I}.$ \qed

\section{A dynamic feedback for the synchronization of trajectories}\label{sec:dynamic}

\noindent In practice, since the switching is not instantaneous
and the tracking error is not perfectly calculated, we have a
serious drawback by using the previous approach, that is the
so-called chattering phenomenon as shown in Example
\ref{exa:static} of Section \ref{sec:static}. To avoid this
phenomenon we propose in the sequel a different coupling term to
add to \eqref{1} in order to solve our synchronization problem,
namely a dynamic feedback coupling term defined by means of the
singular perturbation theory. This method to eliminate the
chattering and to make a controlled dynamical system insensitive
with respect to external perturbations has been introduced in
\cite{CDN1}, where applications to specific tracking problems have
been also presented. Specifically, given the reference trajectory
$y_0(t),$ $t\ge 0,$ which is a solution of system \eqref{2}, we
introduce the function
$s:\mathbb{R}_+\times\mathbb{R}^n\times\mathbb{R}^n\to\mathbb{R}^n$
as follows
$$
  s(t,\xi_0,x)=e^{Ct}(\xi_0-y_0(0))-(x-y_0(t)),
$$
where $C$ is a given real symmetric matrix such that ${\rm
max}\,\lambda(C)\le-\alpha<0,$ here $\lambda(C)$ denotes the set
of the eigenvalues of the matrix $C.$ Then, under conditions
$(3.1)-(3.3),$ for $\varepsilon>0$ small we consider the system
\begin{equation}\label{5}
\left\{
\begin{array}{rcl}
\dot x  &=& \phi(t,x)-Bu\\
\varepsilon\dot u  &=& \dfrac{\partial s}{\partial
t}(t,\xi_0,x)+\dfrac{\partial s}{\partial
x}(t,\xi_0,x)[\phi(t,x)-Bu]:=g(t,\xi_0,x,u).
\end{array}\right.
\end{equation}

\noindent We are now in the position to prove the following
result.
\begin{thm}\label{thm:dynamic}
Assume, that $B$ is a real
negative defined $n\times n$ matrix. Then for any $\delta>0$ there
exists $\varepsilon>0$ such that
$$
  \lim_{t\to\infty}|x_\varepsilon(t)-y_0(t)|<\delta,
$$
where $x_\varepsilon(t), t\ge 0,$ is the solution of (\ref{5})
such that $x_\varepsilon(0)=\xi_0.$
\end{thm}
{\bf Proof.} Let $\varepsilon=0$ and consider the corresponding system
$$
\left\{\begin{array}{l} \dot x=\phi(t,x)-Bu \\
0=Ce^{Ct}(\xi_0-y_0(0))+\dot y_0(t)-\phi(t,x)+Bu,
\end{array}\right.
$$
Resolving the second algebraic equation with respect to $u$ and
substituting in the first equation we obtain
$$
 \dot x_0(t)-\dot y_0(t)=Ce^{Ct}(\xi_0-y_0(0)),\quad t\ge 0
$$
or equivalently
$$
 x_0(t)-y_0(t)=e^{Ct}(\xi_0-y_0(0)),\quad t\ge 0,
$$
with $x_0(t)$ satisfying $x_0(0)=\xi_0.$ By our assumption on the
matrix $C$ we have that $\lim_{t\to\infty}(x_0(t)-y_0(t))=0.$ Now
we show that for given $\delta>0$ there exists $\varepsilon>0$
such that
$$
  |x_\varepsilon(t)-x_0(t)|\le\delta
$$
for any $t\ge 0.$ For this, we use the singular perturbation
theory on unbounded time intervals. Specifically, in the sequel we
verify that all the assumptions of the Theorem at p. 523 of
\cite{hop1} are satisfied for system (\ref{5}) in order to
conclude that
$$
  \lim_{\varepsilon\to 0}x_\varepsilon(t)=x_0(t)
$$
uniformly in $[0,+\infty).$ First of all, for any $(\hat t,\hat
x)\in \mathbb{R}_+\times \mathbb{R}^n$ the equilibrium point
$$
 u(\hat t,\hat x)=-B^{-1}\left[Ce^{Ct}(\xi_0-y_0(0))+
 \psi(\hat t,y_0(\hat t))-\phi(\hat t,\hat x)\right]
$$
of the equation $\varepsilon du/d\tau=g(\hat t,\xi_0,\hat x,u)$ is
asymptotically stable. In other words, the solution $z=z(\tau),$
$\tau\ge 0,$ of the Cauchy problem (the boundary layer)
$$
\left\{\begin{array}{lll} \dot z  =  g(\hat t,\xi_0, \hat x,z)\\
z(0)= z_0. \end{array}\right.
$$
converges asymptotically to $u(\hat t,\hat x).$ In fact, consider
the Lyapunov function $\mathbb{R}_+\times \mathbb{R}^n \times
\mathbb{R}^n \to \mathbb{R}_+$ given by
$$
  V(\hat t,\hat x,z)=\frac{|z-u(\hat t,\hat x)|^2}{2}.
$$
We have, that
$$
\begin{array}{rcl}
  \dfrac{d}{d\tau} V(\hat t,\hat x,z(\tau)) &= &
  \left<\dfrac{\partial}{\partial z}V(\hat t,\hat x,z(\tau)),g(\hat t,\xi_0,\hat x,z(\tau))\right>\\
  &=&\left<z(\tau)-u(\hat t,\hat x),g(\hat t,\xi_0,\hat x,z(\tau))-
  g(\hat t,\xi_0,\hat x,u(\hat t,\hat x))\right>\\
  &= &\left<z(\tau)-u(\hat t,\hat x), B(z(\tau)-u(\hat t,\hat x))\right>\\
  &\le&-\nu|z(\tau)-u(\hat t,\hat x)|^2
\end{array}
$$
for some $\nu>0,$ since $B$ is negative defined. Moreover, observe
that the asymptotic stability is exponential and uniform with
respect to $(\hat t,\hat x)$ when $\hat x$ belongs to bounded
sets.

\noindent Furthermore, the origin $x=0$ is an uniformly stable
equilibrium point of (\ref{5}) when $\varepsilon=0.$ To show this
it is sufficient to observe that for $\varepsilon=0$ we have
$$
  \dot x(t)-\dot y_0(t)=C e^{Ct}(\xi_0-y_0(0)), \qquad t\ge 0
$$
and the change of variable $e=x-y_0$ makes the origin $e=0$ of
this dynamics (exponentially) asymptotically stable, in fact
$$
  \left<e,\dot e\right>=\left<e,C e\right>\le-\alpha|e|^2.
$$
At this point all the assumptions of the Theorem at p. 523 of
\cite{hop1} are satisfied, thus
$$
\lim_{\varepsilon\to 0} x_\varepsilon(t)=x_0(t)
$$
uniformly in $[0,+\infty)$ and the conclusion of the theorem
easily follows.\qed

\begin{rem}\label{rem:2}
The convergence theorem of \cite{hop1} employed in the proof of
Theorem 4 also establishes that the absolutely continuous function
$u_\varepsilon, \varepsilon>0,$ given by \eqref{5} is such that
$$
  \lim_{\varepsilon\to 0} u_\varepsilon(t)=u_0(t)
$$
uniformly on any interval $[t_1,+\infty),$ $t_1>0,$ where
$$
u_0(t)=-B^{-1}\left[Ce^{Ct}(\xi_0-y_0(0))+\psi(t,y_0(t))-\phi(t,x_0(t))\right],
\qquad t\ge 0,
$$
is the so-called equivalent control of the theory of variable
structure systems (see \cite{utkin1} and \cite{utkin2}). In other
words, it is the ideal control which realizes exactly the
condition $s(t,x_0(t))=0.$ This convergence property together with
the absolute continuity of the functions $u_\varepsilon$ and $u_0$
prevents the chattering phenomenon.
\end{rem}

\noindent Finally, we would like to point out that the proposed
dynamical feedback coupling term $-Bu$, designed by means of
\eqref{5}, makes the dynamical system $\dot x=\phi(t,x)$
insensitive with respect to (possibly unknown) perturbations. In
fact, consider any bounded perturbation $p:\mathbb{R}_+ \times
\mathbb{R}^n\to \mathbb{R}^n$ of system $\dot x=\phi(t,x)$
satisfying $(3.1)-(3.2)$. Thus system \eqref{5} takes the form
\begin{equation}\label{6}
\left\{
\begin{array}{rcl}
\dot x & =&\phi(t,x)+p(t,x)-Bu\\
\varepsilon\dot u &=&\dfrac{\partial s}{\partial
t}(t,\xi_0,x)+\dfrac{\partial s}{\partial
x}(t,\xi_0,x)[\phi(t,x)+p(t,x)-Bu].
\end{array}\right.
\end{equation}
It is easy to verify that for this system we can repeat the
arguments used in the proof of Theorem 3. In fact, for
$\varepsilon=0,$ the dynamics of the error $e=x-y_0$ for system
\eqref{6} is the same of that of system \eqref{5}, namely
$$
  \dot x(t)-\dot y_0(t)=C e^{Ct}(\xi_0-y_0(0)), \qquad t\ge 0,
$$
Moreover, $u(\hat t,\hat x)$ is still an asymptotically stable
equilibrium point of the second equation of system \eqref{6}.

\section{Examples}\label{sec:examples}

\noindent In the following examples we have chosen to show the
behavior of only the first component of the solutions of the slave
and master systems, the other components having a quite similar
behavior.

\begin{figure}[h]
\begin{tabular}{l}
\!\!\!\!\!\!\!\!\!\!\!\!\!\!\!\!\!\!\!\!\!\!\!\!\!\!
\includegraphics[scale=.51]{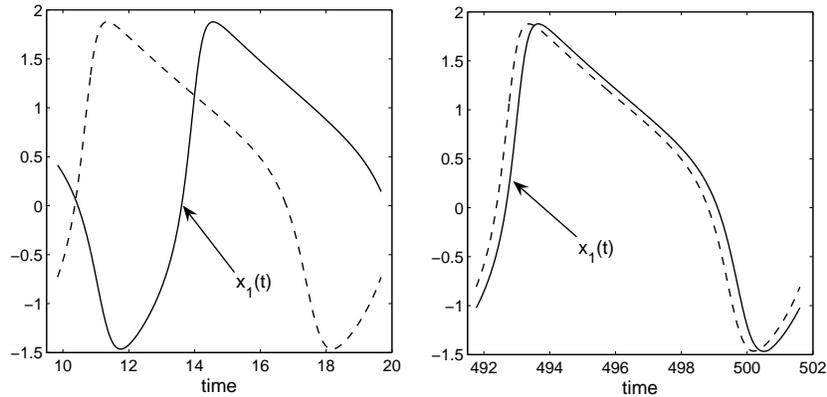}
\end{tabular}
\caption{\footnotesize Time behavior of $x_1(t)$ (solid line) in
Example \ref{exa:phase} in the intervals $[T,2T]$ (left picture)
and $[50T,51T]$ (right picture). The dashed line is the graph of
the first component of the $T$-periodic solution $y_0(t)$ of the
master system.} \label{fig:phase}
\end{figure}

\begin{exa}\label{exa:phase}
At first we illustrate the result of Section \ref{sec:phase} on
the phase synchronization in $\mathbb{R}^2$. We take the same
system as slave and master systems, namely
\begin{equation}\label{eq:phase}
\begin{cases}
\dot{x}_1 = 2\left(x_1-\dfrac{1}{3}x_1^3+x_2-\dfrac{9}{20}\right)\vspace{4pt}\\
\dot{x}_2 = -\dfrac{1}{2}\left(x_1+\dfrac{4}{5}x_2-\dfrac{7}{10}\right)
\end{cases}
\end{equation}
which has an asymptotically stable limit cycle with period $T
\simeq 9.83$. In particular, let $y_0$ be the $T$-periodic
solution of \eqref{eq:phase} which starts at $t=0$ from the point
$(-0.7481,1.5164)$ of the limit cycle. In this way, since
$$
\min_{s\in[0,T]}\frac{1}{T}\int_0^T | y_0(\tau+s)-y_0(\tau)|^2
d\tau=0,\quad \mbox{system \eqref{peq} becomes}
$$
\begin{equation}\label{eq:phaseslave}
\begin{cases}
\dot{x}_1 = 2\left(x_1-\dfrac{1}{3}x_1^3+x_2-\dfrac{9}{20}\right)
(1+\varepsilon(|x-y_0(t)|^2-\delta)
\vspace{4pt}\\
\dot{x}_2 = -\dfrac{1}{2}\left(x_1+\dfrac{4}{5}x_2-\dfrac{7}{10}\right)
(1+\varepsilon(|x-y_0(t)|^2-\delta)
\end{cases}
\end{equation}
and we consider the solution $x$ of \eqref{eq:phaseslave} which starts from $x(0)=(5,-5),$ with the choice
$\varepsilon=0.01$ and $\delta=0.05.$
In Figure \ref{fig:phase} the solid line is the graph of the first component $x_1$ of $x$ in the intervals $[T,2T]$
(left picture) and $[50T,51T]$ (right picture), while the dashed line is the graph of the first component of the
$T$-periodic solution $y_0.$
Due to the asymptotical stability of the limit cycle of \eqref{eq:phase}, we see that $x_1$ already has
the shape of the first component of $y_0$ after only one period, but the phases of these trajectories are still
considerably different.
The picture on the right shows that the phase difference between the first components of $x$ and $y_0$ is
significantly reduced after $50$ periods thanks to the coupling term added in \eqref{eq:phaseslave}.
\end{exa}

\begin{figure}
\begin{tabular}{l}
\!\!\!\!\!\!\!\!\!\!\!\!\!\!\!\!\!\!\!\!\!\!\!
\includegraphics[scale=.51]{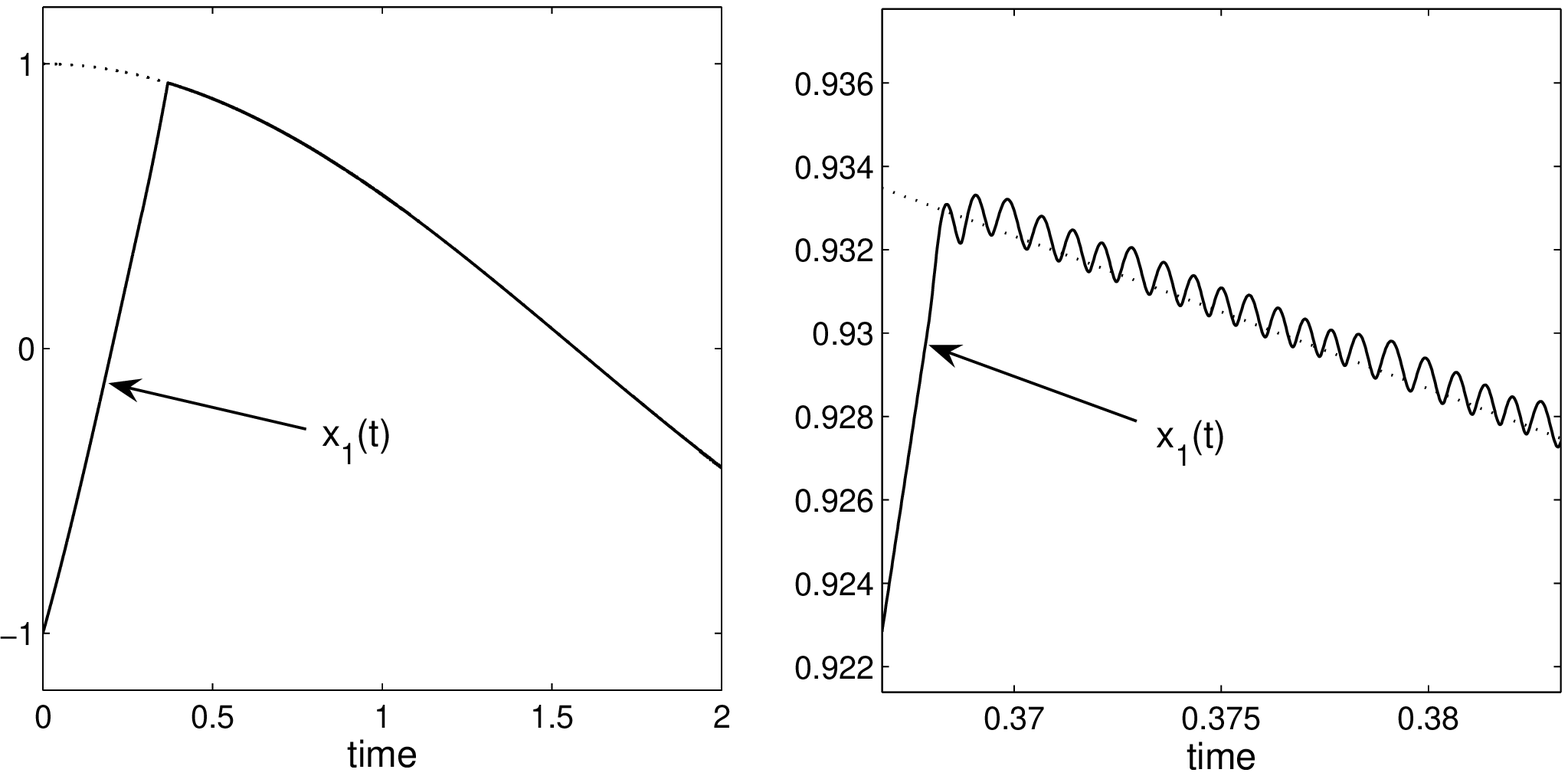}
\end{tabular}
\caption{\footnotesize Time behavior of $x_1(t)$ (solid line) of
Example \ref{exa:static}. The dotted line is the graph of $\cos
t,$ that is the first component of the periodic solution of the
master system \eqref{eq:master}.} \label{fig:static}
\end{figure}

\begin{exa}\label{exa:static}
We exploit the result of Section \ref{sec:static} by considering the static feedback synchronization of
the following neural network introduced in \cite{ZN}
\begin{equation}\label{eq:slave}
\dot{x} =
-\left[\begin{array}{ccc}1&0&0\\0&1&0\\0&0&1\end{array}\right]x+
\left[\begin{array}{lll}
\hphantom{-}1.25&-3.2&-3.2\\
-3.2&\hphantom{-}1.1&-4.4\\
-3.2&\hphantom{-}4.4&\hphantom{-}1\end{array}\right]
\left[\begin{array}{l} f(x_1)\\f(x_2)\\f(x_3)\end{array}\right],
\end{equation}
where $f(s)=(|s+1|-|s-1|)/2.$
In \cite{ZN} it is shown that system \eqref{eq:slave} has a chaotic attractor.
As master system we use another neural network
\begin{equation}\label{eq:master}
\dot{x} =
-\left[\begin{array}{ccc}\frac{10}{7}&0&0\\0&1&0\\0&0&0.1\end{array}\right]x+
\left[\begin{array}{rrc}
-\frac{20}{7}&10&0\\
1&-30&1\\
0&\frac{100}{7}&-1.9\end{array}\right]
\left[\begin{array}{l} x_1^3\\x_2\\x_3\end{array}\right] + I(t),
\end{equation}
where $I(t)$ is the periodic input such that \eqref{eq:master} has
the $2\pi$-periodic solution $y_0(t)=(\cos t,\sin t,-\cos t).$ We
take the gains $b_1=b_2=b_3=-3.5$ in the coupling term \eqref{3}
and consider the solution $x$ of \eqref{1p} such that
$x(0)=(-1,1,1);$ then in Figure \ref{fig:static} the graph of the
first component $x_1$ of $x$ is plotted. In particular, the
picture on the left shows the convergence in finite time to the
first component of the reference solution $y_0$; on the other
hand, the picture on the right is a zoom of the left one around
the hitting zone and is a clear evidence of the chattering
phenomenon.
\end{exa}

\begin{figure}
\begin{tabular}{l}
\!\!\!\!\!\!\!\!\!\!\!\!\!\!\!\!\!\!\!\!\!\!\!
\includegraphics[scale=.51]{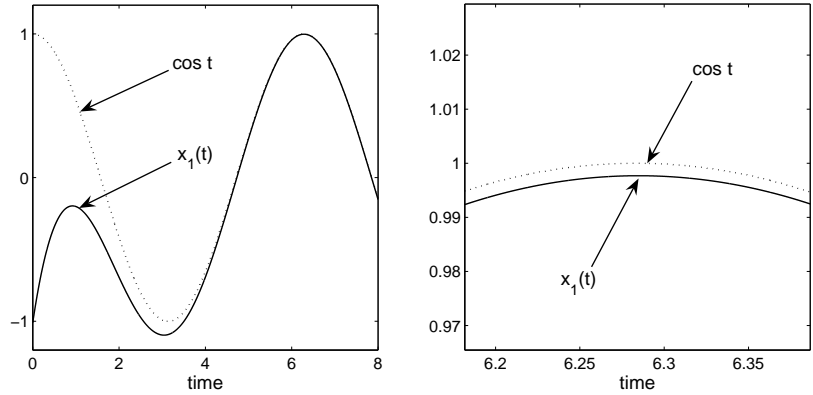}
\end{tabular}
\caption{\footnotesize Time behavior of $x_1(t)$ (solid line) of
Example \ref{exa:dynamic}. The dotted line is the graph of $\cos
t,$ that is the first component of the periodic solution of the
master system.}
  \label{fig:dynamic}
\end{figure}

\begin{exa}\label{exa:dynamic}
To illustrate the dynamic feedback synchronization considered in
Section \ref{sec:dynamic} and to compare the issues with those of
the previous example, we consider again the two neural networks
\eqref{eq:slave} and \eqref{eq:master} of Example \ref{exa:static}
with the same forcing term $I(t).$ We take $C=B=\diag(-1,-1,-1)$
and $\varepsilon=0.001$ and consider the solution $x$ of \eqref{5}
such that $x(0)=(-1,1,1).$ In Figure \ref{fig:dynamic} the graphs
of the first components $x_1$ of $x$ (solid line) and $\cos t$ of
$y_0$ (dotted line) are plotted. In particular, the picture on the
left shows how $x_1$ approaches $\cos t,$ while the picture on the
right is a zoom of the previous one around $t=2\pi$ and shows that
$x_1(t)$ remains in a neighborhood of $\cos t$ without any
chattering.
\end{exa}

\begin{rem}\label{rem:3}
Observe that in Examples \ref{exa:static} and \ref{exa:dynamic} we
have taken the initial conditions for the master system on the
reference trajectory $y_0.$ On the other hand, in \cite{PT} it is
shown that $y_0$ is globally exponentially stable, therefore these
initial conditions could be also chosen as far as we like from
$y_0.$
\end{rem}

\section*{Acknowledgment}
The first author was supported by a President of Russian
Federation Fellowship for Scientific Training Abroad and by the
Grant A04-2.8-64 of Russian Federation Federal Agency on
Education. The second and third author were supported by the
research project GNAMPA :``Qualitative Analysis and Control of
Hybrid Systems''.

\end{document}